\documentclass[11pt]{amsart}
\usepackage{amscd}      % to be able to use "CD" environment
\usepackage{amssymb}
\usepackage{amsmath}
\usepackage{xypic}      % to be able to use "diagram" environment of
\LaTeXdiagrams          % xypic to typeset commutative diagrams
\usepackage[all,v2]{xy}
\xyoption{2cell} \UseAllTwocells \xyoption{frame} \CompileMatrices
\allowdisplaybreaks[3]

\usepackage{latexsym}
\usepackage{epsfig}
\usepackage{amsfonts}
\usepackage{enumerate}

\newtheorem{prop}{Proposition}[section]
\newtheorem{lem}[prop]{Lemma}

\newtheorem{cor}[prop]{Corollary}
\newtheorem{thm}[prop]{Theorem}
\newtheorem{rmk}[prop]{Remark}
\newtheorem{example}[prop]{Example}

\newtheorem{defn}[prop]{Definition}

\def \jeden {1\hskip-3.5pt1}

\newenvironment{pf}{\begin{trivlist}\item[]{\sc Proof.}}%
            {\nolinebreak $\Box$ \end{trivlist}}

\newcommand{\noprint}[1]{}

\renewcommand{\tilde}{\widetilde}

\newcommand{\YY}{{\mathfrak Y}}

\newcommand{\FF}{{\mathfrak F}}

\newcommand{\MM}{{\mathfrak M}}

\newcommand{\Cc}{{\mathfrak c}}

\newcommand{\zz}{{\mathbb Z}}

\newcommand{\aaa}{{\mathbb A}}

\newcommand{\pp}{{\mathbb P}}
\newcommand{\cc}{{\mathbb C}}

\renewcommand{\O}{{\mathcal O}}

\DeclareMathOperator{\dg}{dg}
\DeclareMathOperator{\val}{val}

\newcommand{\ord}{\mathop{\rm ord}\nolimits}

\def\Label#1{\label{#1}{\tt [#1]}\phantom{h}}
\def\Label{\label}

\title{Motivic Milnor fibre of cyclic $L_\infty$-algebras}
\author{Yunfeng Jiang}
\address{Department of Mathematics\\ University of Utah\\ 155 S 1400 E JWB 233\\Salt Lake city\\ UT 84112\\ USA} 
\email{jiangyf@math.utah.edu}

\begin{document}
\sloppy \maketitle
\begin{abstract}
We define motivic Milnor fiber of cyclic $L_\infty$-algebras of dimension three using the method of 
Denef and Loeser of motivic integration.  
It is proved by Nicaise and Sebag that the topological Euler characteristic of the motivic Milnor fiber is equal to the 
Euler characteristic of the \'{e}tale cohomology of  the analytic Milnor fiber. 
We prove that the value of Behrend function on the germ  moduli space determined by $L$ is equal to the Euler characteristic of the analytic Milnor fiber.  Thus we prove that  Behrend function depends only on the formal neighborhood of  the moduli space.
\end{abstract}

\tableofcontents

%%% -----------------------------------------------------------------------
\maketitle
%%% ----------------------------------------------------------------------

\section{Introduction}

Let $L$ be an $L_\infty$-algebra with finite dimensional cohomology.  If there is a 
nondegenerate
symmetric 
bilinear form on $L$ which has degree $3$, see definition in Section \ref{l-infinity}, we call 
$L$ a cyclic $L_\infty$-algebra. From transfer theorem in \cite{KS}, \cite{BG}, there is  a 
cyclic $L_\infty$-algebra structure on the cohomology $H(L)$.  In this paper we assume that 
$H^{i}(L)=0$ only except for $i=1,2$. Then by cyclic property $H^{1}(L)\cong H^{2}(L)^{\vee}$.

From the cyclic $L_\infty$-algebra structure on $H(L)$,  there is a potential function 
$f: H^{1}(L)\rightarrow \cc$, which is a formal power series beginning from cubic 
terms.  Let $M:=H^{1}(L)$ and let $0\in M$ be the origin. The germ (formal) moduli space 
$(X,0)$ associated to $L$ is the critical locus of $f$, i.e. the vanishing locus of $df$.

Following \cite{DeLo1}, \cite{KS2}, we define motivic Milnor fiber $MF_0(L)$ of $L$ using 
resolution of singularities of the critical subscheme of $df$.  This motivic  Milnor fiber
$MF_0(L)$ can be seen as an element in the equivariant Grothendieck group of varieties 
over $\cc$. 
Let $X_0$ be the zero locus of $f$, and for each 
$x\in X_0$, we define the  motivic Milnor fiber $MF_x(L)$ at $x$ in a similar way.
Taking Euler characteristic we get a constructible function 
on $X_0$.

We also define  analytic Milnor fiber of $L$ in the sense of \cite{NS}. This analytic Milnor fiber is an analytic space over the non-archimedean field $K=\cc((t))$.  From \cite{NS}, the Euler characteristic of the \'{e}tale cohomology of the analytic Minor fiber in the sense of Berkovich \cite{Ber}, is equal to the topological Euler characteristic of the motivic Milnor fiber. Using this result we prove that the value of   Behrend function on the critical set $X=\mathbb{V}(df)$ is given by  the  Euler characteristic of the analytic Minor fiber.
This way we prove that the Behrend function depends only on  the formal neighborhood of the points in the moduli space.

Moreover, if $L$ is a Donaldson-Thomas type dg Lie algebra, then the motivic Milnor fiber can give 
the pointed Donaldson-Thomas invariant determined by $L$.    In the paper \cite{BG}, Behrend and Getzler prove that the formal potential function $f$ for the cyclic $\dg$ Lie algebra $L$ coming from the Schur objects in the derived category of coherent sheaves on Calabi-Yau threefolds can be made to be convergent over a local neighborhood of the origin $0\in M$.
In this case, the pointed Donaldson-Thomas invariant determined by $L$ is given by the general Milnor number of $f$ in \cite{Mil}, see \cite{Be}.

The paper is outlined as follows.  We briefly review the definition of cyclic $L_\infty$-algebra and deformation theory in Section \ref{l-infinity}. For more details about cyclic $L_\infty$-algebras, see \cite{BG}. In Section \ref{MC},  the Maurer-Cartan space of a cyclic $L_\infty$-algebra is discussed, and the Maurer-Cartan space gives the germ moduli space of $L$. We define motivic Milnor fiber of cyclic $L_\infty$-algebras in Section \ref{motivic-milnor}. Finally in Section \ref{behrend}, we relate the Behrend function of the moduli space of $L$ to motivic Milnor fiber of $L$.

\subsection*{Convention}
 
 We work over $\cc$ throughout the paper. Denote by $R=\cc[[t]]$ the ring of formal power series, and 
 $K=\cc((t))$, the quotient field. For the analytic Milnor fiber $\FF_{x}$, which is an analytic space  over the non-archimedean field $K$, we use the \'{e}tale cohomology of  $\FF_{x}$ in the sense of Berkovich, \cite{Ber}.

%%% ----------------------------------------------------------------------
\subsection*{Acknowledgments}

The author would like to thank Professors Kai Behrend and Yan Soibelman  for their encouragements and answering  related questions.  The author is very grateful to Professor Johannes Nicaise for 
giving valuable suggestions on analytic Milnor fibers.

%%%%%%%%%%%%%%%%%%%%%%%%%%%%%%%%%%%%%%%%%%%%%%%%%%%%%%%%%%%%%%%%%
%%% ----------------------------------------------------------------------
\section{The cyclic $L_\infty$-algebra and deformation theory.}\Label{l-infinity}
%%%%--------------------------------------------------------------

\subsection{Cyclic $L_{\infty}$-algebras.}

\begin{defn}\label{l-infty}
An $L_{\infty}$-algebra $L$ is a $\mathbb{Z}$-graded $\cc$-vector space
$\oplus_{i}L^{i}$ equipped with linear maps:
$$\mu_{k}: \Lambda^{k}L\longrightarrow L[2-k],$$
given by 
$a_{1}\otimes\cdots\otimes a_{k}\longmapsto \mu_{k}(a_{1},\cdots,a_{k})$ for $k\geq 1$,
which satisfies the $higher~order~Jacobi~identities$ for any $a_{1},\cdots,a_{n}$:
\begin{align}\label{jacobi-linfty}
&\sum_{l=1}^{n}\sum_{\sigma\in Sh(l,n-l)}(-1)^{\widetilde{\sigma}+(n-l+1)(l-1)}
\epsilon(\sigma;a_{1},\cdots,a_{n})\cdot \nonumber \\
&\mu_{n-l+1}(\mu_{l}(a_{\sigma(1)},\cdots,a_{\sigma(l)}),a_{\sigma(l+1)},\cdots,a_{\sigma(n)})=0,
\end{align}
where the shuffle $Sh(l,n-l)$ is the set of all permutations $\sigma: \{1,\cdots,n\}\rightarrow\{1,\cdots,n\}$
satisfying  $\sigma(1)<\cdots<\sigma(l)$ and  $\sigma(l+1)<\cdots<\sigma(n)$.
The symbol $\epsilon(\sigma;a_{1},\cdots,a_{n})$ (which we abbreviate $\epsilon(\sigma)$) stands for the 
$Koszul~sign$ defined by
$$a_{\sigma(1)}\wedge\cdots\wedge a_{\sigma(n)}=(-1)^{\widetilde{\sigma}}\epsilon(\sigma)a_{1}\wedge\cdots\wedge a_{n},$$
where $\widetilde{\sigma}$ is the parity of the permutation $\sigma$.
\end{defn}

%%%--------------------------------------------------------
%%%------------------------------
\textbf{The $L_{\infty}$-morphisms.}

\begin{defn}\label{linfty-morphism}
An $L_{\infty}$ morphism between two $L_{\infty}$-algebras $L$ and $L^{'}$
is defined by a degree zero coalgebra morphism $F$ from $S(L[1])$ to
$S(L^{'}[1])$ which commutes with the codifferentials $Q$ and $Q^{'}$.
It is completely determined by a set of linear maps $F_{n}: S^{n}(L[1])\rightarrow L^{'}[1-n]$
(or equivalently $F_{n}: \bigwedge^{n}L\rightarrow L^{'}[1-n]$) satisfying a 
set of equations.
\end{defn}

An $L_{\infty}$-morphism $F: L\longrightarrow L^{'}$ is called quasi-isomorphic if the first
component $F_{1}: L\longrightarrow L^{'}$ induces an isomorphism between cohomology 
groups of complexes $(L,\mu_{1})$ and $(L^{'},\mu_{1}^{'})$.

%%%------------------------------------
\textbf{The cyclic $L_{\infty}$-algebras.}

\begin{defn}\label{cyclic-l-infty}
A  \textbf{Cyclic $L_{\infty}$-algebra} of dimension $d$ is a triple  $(L,\mu,\ae)$, where 
$(L,\mu)$ is an $L_{\infty}$-algebra with linear maps
$(\mu_{k})$, and 
$$\ae: L\otimes L\longrightarrow \mathbb{C}[-d]$$
is a perfect pairing which means that 
$$H^{i}(L)\otimes H^{d-i}(L)\longrightarrow \mathbb{C}$$
is a perfect pairing of finite dimensional vector spaces  for each $i$. The 
bilinear form $\ae$ satisfies the following conditions:
\begin{enumerate}
\item $\ae$ is graded symmetric, i.e. $\ae(a,b)=(-1)^{ab}\ae(b,a)$;
\item For any $n\geq 1$, 
$$\ae(\mu_{n}(x_1,\cdots,x_n),x_{n+1})=(-1)^{n+x_1(x_2+\cdots+x_{n+1})}\ae(\mu_{n}(x_2,\cdots,x_{n+1}),x_1).$$
\end{enumerate}
\end{defn}

%%%---------------------------------------------------------------------------
\textbf{Transfer theorem.}

Let $(L,\mu,\ae)$ be a cyclic $L_{\infty}$-algebra.  We write $d=\mu_1$. Let
$$\eta: L\longrightarrow L[-1]$$
be map of degree $-1$, such that 
\begin{enumerate}
\item $\eta^{2}=0$;
\item $\eta d\eta=\eta$;
\item $\ae(\eta x,y)+(-1)^{x}\ae(x,\eta y)=0$.
\end{enumerate}
Let $$\Pi=1-[d,\eta],$$
where $[d,\eta]=d\eta+\eta d$. Then $\Pi^{2}=\Pi$. 

Let
$$M=\Pi(L).$$
Let 
$\iota: M\rightarrow L$ be the inclusion and $p:L\rightarrow M$ the projection. Then
$$\iota p=\Pi, ~~ p\iota=id_{M}.$$

\begin{thm}\label{linfty-coh}(\cite{BG})
Let $(L,\mu,\ae)$ be a cyclic $L_{\infty}$-algebra, then there exists a cyclic  $L_{\infty}$-algebra structure on 
$M$ such that there is a $L_{\infty}$-morphism 
$$\varphi: M\longrightarrow L$$
as $L_{\infty}$-algebras.
\end{thm}

\begin{rmk}
For the cohomology $H^\bullet(L)$ of $L$, one can choose the contraction operator $\eta$
such that $\varphi:  H^\bullet(L) \longrightarrow L$ is a quasi-isomorphic morphism of $L_\infty$-algebras.
\end{rmk}

%%%-------------------------------------------------------
\textbf{Cyclic differential graded Lie algebras}\\
A special case of cyclic $L_\infty$-algebra is the cyclic $\dg$ Lie algebra.
\begin{defn}\label{dgla}
A \textbf{Differential Graded Lie Algebra} ($\dg$ Lie algebra) is a triple $(L,[,],d)$,
where $L=\oplus_{i}L_{i}$ is a $\mathbb{Z}$-graded vector space,  and the bracket
$[,]: L\times L\longrightarrow L$ is  bilinear and 
$d: L\longrightarrow L$ is a linear map such that 
\begin{enumerate}
\item The bracket $[,]$ is homogeneous skewsymmetric which means that 
$[L^{i},L^{j}]\subset L^{i+j}$ and $[a,b]+(-1)^{\overline{a}\overline{b}}[b,a]=0$
for every $a, b$ homogeneous.
\item Every homogeneous elements $a, b, c$ satisfy the Jacobi identity:
$$[a,[b,c]]=[[a,b],c]+(-1)^{\overline{a}\overline{b}}[b,[a,c]].$$
\item The map $d$ has degree 1,  $d^{2}=0$ and $d[a,b]=[da,b]+(-1)^{\overline{a}}[a,db]$,
$d$ is called the differential of $L$.
\end{enumerate}
\end{defn}

\begin{defn}\label{cyclicdgla}
A \textbf{Cyclic Differential Graded Lie Algebra of dimension $d$} is a DGLA
$L$ together with a bilinear form
$$\ae: L\times L\longrightarrow \cc[d]$$
such that for any homogeneous elements $a, b, c$,
\begin{enumerate}
\item $\ae$ is graded symmetric, i.e. 
$\ae(a,b)=(-1)^{\overline{a}\overline{b}}\ae(b,a)$;
\item $\ae(da,b)+(-1)^{\overline{a}}\ae(a,db)=0$;
\item $\ae([a,b],c)=\ae(a,[b,c])$;
\item $\ae$ induces a perfect pairing
$$H^{d}(L\otimes L)\longrightarrow \mathbb{C}$$
which induces a perfect pairing 
$$H^{i}(L)\otimes H^{d-i}(L)\longrightarrow \mathbb{C}$$
for any $i$.
\end{enumerate}
\end{defn}

For a cyclic $\dg$ Lie algebra $L$, let
$$Z^{i}=ker(d^{i}: L^{i}\rightarrow L^{i+1})$$
and
$$B^{i}=im(d^{i-1}: L^{i-1}\rightarrow L^{i}),$$ 
then the 
cohomology is $H^{i}(L)=Z^{i}/B^{i}$. The cohomology 
$H^{\bullet}=\oplus_{i}H^{i}(L)$ has a induced $\dg$ Lie algebra structure
with zero differential. 

\begin{defn}
A morphism $f: L\longrightarrow L^{'}$ betrween two $\dg$ Lie algebras
is a morphism as graded vector spaces compatible with the 
differential and the Lie bracket. The morphism $f$ is called a
quasi-isomorphism if it induces an isomorphism on the cohomology $\dg$ Lie algebras. 
\end{defn}
\begin{example}
Let $Y$ be a smooth Calabi-Yau 3-fold. Let $\mathcal{F}$ be any derived object in the derived category of coherent sheaves on $Y$. We take $E^{\bullet}\rightarrow \mathcal{F}$ as a projective resolution. Let 
$$L=\bigoplus_n\bigoplus_{p+q=n}\mathcal{A}^{0,p}(X,\mathcal{H}om^q(E^\bullet,E^\bullet)_0),$$
where $\mathcal{H}om^q(E^\bullet,E^\bullet)_0$ is taken as sheaf of vector bundles over $X$. On the complex
$\mathcal{H}om^q(E^\bullet,E^\bullet)_0$ there exists a bracket given by
$$[\varphi,\psi]=\varphi\psi-(-1)^{\varphi\psi}[\psi,\varphi].$$
We define the differential 
$$d=\overline{\partial},$$
on $L$, where $\overline{\partial}$ is the Dolbeault differential. The bracket $[,]$ is defined by
$$[\omega\otimes\varphi,\tau\otimes\psi]=\omega\wedge \tau\otimes [\varphi,\psi].$$
Then we check that $L$ is a $\dg$ Lie algebra under $d$ and $[,]$.

Again using the trace on $\mathcal{H}om(E^\bullet,E^\bullet)$, let
$$\ae: \mathcal{H}om(E^\bullet,E^\bullet)\times \mathcal{H}om(E^\bullet,E^\bullet)\longrightarrow \mathcal{O}_X$$
be given by
$$\varphi\otimes\psi\longmapsto \mbox{tr}(\varphi\circ \psi).$$
Then $\ae$ is symmetric and induces
$$\ae: L\times L\longrightarrow \mathbb{C}[-3]$$
by
$$\mathcal{A}^{0,j}(X,\mathcal{H}om^q(E^\bullet,E^\bullet)_0)\otimes \mathcal{A}^{0,3-j}(X,\mathcal{H}om^q(E^\bullet,E^\bullet)_0)
\longrightarrow \mathbb{C}$$
and
$$\varphi\otimes \psi\longmapsto \int_{X}\mbox{tr}(\varphi\circ \psi)\wedge \omega_X.$$
We check that $\ae$ is symmetric and satisfies the condition in Definition \ref{cyclicdgla}.
So $L$ is a cyclic $\dg$ Lie algebra of dimension $3$, which we call Donaldson-Thomas type $\dg$
Lie algebra.
\end{example}

%%%-----------------------------------------------------------------------------------------------------
\subsection{Deformation theory of $L_{\infty}$-algebras.}
For a $\dg$ Lie algebra $L$ the deformation theory of $L$ is given by the Maurer-Cartan equation.
$L_{\infty}$-algebras are generalizations of $\dg$ Lie algebras. We have a modified deformation functor 
for the $L_{\infty}$-algebras.

Let $Art_{\cc}$ be the category of local Artin algebras over $\cc$.
\begin{defn}
The modified deformation functor $Def_{L}: Art_{\cc}\longrightarrow Sets$ of an $L_{\infty}$-algebra 
is given by:
$$Def_{L}(B)=\{a\in C(L^{1}\otimes m_{B})|\sum_{i=1}^{\infty}\tfrac{(-1)^{i(i+1)/2}}{i!}
\mu_{i}(a^{\otimes i})=0\},$$
where  $B\in Art_{\cc}$, $C(L^{1}\otimes m_{B})$ is the complement of $d(L^{0}\otimes m_{B})$
in $L^{1}\otimes m_{B}$. This is also called the homotopy Maurer-Cartan functor.
\end{defn}

Let $F: L\longrightarrow L^{'}$ be a morphism between $L_{\infty}$-algebras
$L$ and $L^{'}$. Then if  $F$ is a quasi-isoporphism, it reduces an isomorphism
$Def_{L}\cong Def_{L^{'}}$ for the deformation functors, see \cite{Man1}.

%%%-------------------------------------------------------------------------------
\section{The Maurer-Cartan space.}\label{MC}

\subsection{The moduli space.}

Let $L$ be a cyclic $L_\infty$-algebra of dimension $3$.  
Let
\begin{equation}
X:=\lbrace{x\in L^{1}|\sum_{i=1}^{\infty}\tfrac{(-1)^{i(i+1)/2}}{i!}
\mu_{i}(x^{\otimes i})=0\rbrace}.
\end{equation}
The space $X$ is called the Maurer-Cartan space associated to the cyclic
$L_\infty$-algebra $L$.

By the transfer theorem there is a cyclic $L_\infty$-algebra structure on the cohomology
$H(L)$.  Let $(H(L),\nu, \kappa)$ be the cyclic $L_\infty$-algebra structure. 
The Maurer-Cartan map
\begin{equation}\label{HMC}
MC=\bigoplus_{k}F_k: \Lambda^{\bullet}H^{1}(L)\longrightarrow H^{2}(L)
\end{equation}
is given by
$$x_{1}\otimes\cdots\otimes x_{k}\longmapsto \sum_{k}\tfrac{(-1)^{k(k+1)/2}}{k!}\nu_{k}(x_{1}\otimes\cdots\otimes a_{k}).$$

The Maurer-Cartan space of $H(L)$ is given by
\begin{equation}\label{mc}
MC=\{ x\in H^{1}(L)| \oplus_{k} F_k(x)=0\}.
\end{equation}

\begin{thm}
If $H^{0}(L)=0$, then the Maurer-Cartan space $MC$ is isomorphic to the 
completion $(\hat{X},0)$ of the Maurer-Cartan space $X$ with respect to 
the origin $0\in L^{1}$.
\end{thm}
\begin{pf}
See \cite{Man1}, \cite{Man2}.
\end{pf}

\subsection{The formal potential.}

In this section we talk about a special  case that the cohomology of the cyclic 
$L_\infty$-algebra $L$ satisfies that $H^{i}\neq 0$ only for $i=1,2$.
Then the structure maps  are given by the following
\begin{equation}
\nu_k: \Lambda^{k}H^{1}(L)\longrightarrow H^{2}(L)
\end{equation}
And the bilinear form $\kappa$ is given by 
$$\kappa: H^{1}(L)\longrightarrow H^{2}(L)=H^{1}(L)^{\vee}.$$

Define the potential  function
\begin{equation}\label{potential}
f: H^{1}(L)\longrightarrow \mathbb{C}
\end{equation}
by
$$f(z)=\sum_{n=2}^{\infty}\frac{(-1)^{n(n+1)/2}}{(n+1)!}\kappa(\nu_{n}(z,\cdots,z),z),$$
for $z\in H^1(L)$.

\begin{prop}(\cite{BG})
We have that $df=F$.
\end{prop}
\begin{pf}
This result is proved in \cite{BG}. For the convenience, we recall the proof here.
Recall that $$F=\oplus_{k}F_k,$$
where $F_k: H^{1}(L)\longrightarrow H^2(L)$ is defined by 
$$F_{k}(z)=\frac{(-1)^{k(k+1)/2}}{k!}\nu_{k}(z,\cdots,z).$$
Since $H(L)$ is a cyclic $L_{\infty}$-algebra, the bilinear form 
$$\kappa: H(L)\times H(L)\longrightarrow H(L)$$
induces an isomorphism
$$\kappa: H^{1}(L)\stackrel{\cong}{\longrightarrow}(H^2(L))^{\vee}.$$
Thus $F_k$ is a defferential form on $H^{1}(L)$ for each $k\geq 2$.
Since
$$F_k(z)\in H^{2}(L)\cong (H^{1}(L))^{\vee},$$
it can be taken as a function on $H^1(L)$, and for $x\in H^{1}(L)$,
$$F_k(z)(x)=\frac{(-1)^{k(k+1)/2}}{k!}\kappa(x,\nu_{k}(z,\cdots,z)).$$

Consider the formal function $f$ defined in (\ref{potential}).
Let $$f=\oplus_{k}f_k,$$
where $f_k$ is the function defined by 
$$f(z)=\frac{(-1)^{k(k+1)/2}}{(k+1)!}\kappa(\nu_{k}(z,\cdots,z),z).$$
Then we have
that for any $x\in H^1(L)$,

\begin{eqnarray}
\frac{df_k}{dz}(x)&=&\sum_{i=1}^{k+1}\frac{(-1)^{k(k+1)/2}}{(k+1)!}\kappa(\nu_{k}(z,\cdots,x,\cdots,z),z)\nonumber \\
&=&(k+1)\frac{(-1)^{k(k+1)/2}}{(k+1)!}\kappa(x,\nu_{k}(z,\cdots,z))\nonumber \\
&=&\frac{(-1)^{k(k+1)/2}}{k!}\kappa(x,\nu_{k}(z,\cdots,z)), \nonumber
\end{eqnarray}
where the second equality comes from the cyclic property of
$\kappa$ and $(\nu_k)_{k\geq 2}$. So 
$$\frac{df_k}{dz}=F_k(z), ~~\mbox{and}~ df=F.$$
\end{pf}

Let $\{x_1,\cdots,x_m\}$ be a basis of $H^{1}(L)^{\vee}$. Let 
$$A=\frac{\mathbb{C}\langle x_1,\cdots,x_m\rangle }{\partial_{e}(f)},$$
where $\partial_{e}(f)$ represents the derivative of the formal potential function $f$.
Then we have the following result:

\begin{thm}
The deformation functor $Def_{L}$ of the cyclic $L_\infty$-algebra $L$ is 
pro-represented by the local complete algebra $A$.
\end{thm}

%%%%%%%%%%%%%%%%%%%%%%%%%%%%%%%%%%%%%%%%%%%%%%

\section{The motivic Milnor fiber.}\label{motivic-milnor}

\subsection{Resolution of singularities.}

Let $(L,\mu, \ae)$ be a cyclic $L_\infty$-algebra of dimension $3$.
By transfer theorem there is a cyclic $L_\infty$-algebra structure on the cohomology
$H(L)$, which we denote by $(H(L), \nu, \kappa)$.

The cyclicity property of $H(L)$ gives a formal potential function
\begin{equation}
f: H^{1}(L)\longrightarrow \mathbb{C}
\end{equation}
defined by 
$$z\mapsto 
f(z)=\sum_{n=2}^\infty\frac{(-1)^{\frac{n(n+1)}{2}}}{(n+1)!}\kappa\big(\nu_n(z,\ldots,z),z\big).
$$

In general $f$ is a formal series and $f(0)=0$, where $0\in H^{1}(L)$ is the origin. We 
generalize the motivic Milnor fiber of the holomorphic functions to the formal 
series $f$ as suggested in \cite{KS2}.
For simplicity, let $M:=H^{1}(L)$ and let $m=\dim(M)$.
Let $x\in X_0$ be a point such that $f(x)=0$, and 
$f$ is convergent at the point $x$. Let 
$$\widehat{f}_{x}: \MM\longrightarrow \cc[[t]]$$
be the formal completion of $f$ at $x$. Then 
$\MM$ is a special formal $\cc[[t]]$-scheme. 
First we have a result of resolution of singularities.

\begin{thm}\label{reso}(\cite{Tem},\cite{NS})
Let $\MM$ be an  affine generically smooth formal $\cc[[t]]$-scheme. 
Then $\MM$ admits a resolution of singularities by means of formal admissible blow-ups.
$\square$
\end{thm}

Using Theorem \ref{reso},  let 
\begin{equation}
h:  \YY\longrightarrow \MM
\end{equation}
be the  resolution  of singularity of formal  $\cc[[t]]$-scheme $\MM$.

Let $E_i$, $i\in A$, be the set of irreducible components of the exceptional divisors of the resolution. 
For $I\subset A$,  we set 
$$E_{I}:=\bigcap_{i\in I}E_{i}$$
and 
$$E_{I}^{\circ}:=E_{I}\setminus \bigcup_{j\notin I}E_j.$$
Let $m_{i}$ be the multiplicity of the component $E_i$, which means that 
the special fiber of the resolution is 
$$\sum_{i\in A}m_iE_i.$$

\subsection{Grothendieck group of varieties.}

In this section we briefly review the Grothendieck group of varieties. 
Let $S$ be an algebraic variety over $\cc$. Let $Var_{S}$ be the category of 
$S$-varieties.

Let $K_0(Var_{S})$ be the Grothendieck group of $S$-varieties.  By definition $K_0(Var_{S})$ 
is an abelian group with generators given by all the varieties $[X]$'s, where $X\rightarrow S$ are $S$-varieties,  and the relations are $[X]=[Y]$, if $X$ is isomorphic to $Y$, and 
$[X]=[Y]+[X\setminus Y]$ if $Y$ is a Zariski closed subvariety of $X$.
Let $[X],  [Y]\in K_0(Var_{S})$,  and define $[X][Y]=[X\times_{S} Y]$.  Then 
we have a product on $K_0(Var_{S})$. 
Let $\mathbb{L}$ represent the class of $[\mathbb{A}_{\cc}^{1}\times S]$.
Let $\mathcal{M}_{S}=K_0(Var_{S})[\mathbb{L}^{-1}]$
be the ring by inverting the class $\mathbb{L}$ in the ring $K_0(Var_{S})$.

If $S$ is a point $Spec (\cc)$, we write $K_0(Var_{\cc})$ for the Grothendieck ring of $\cc$-varieties.
One can take the map $Var_{\cc}\longrightarrow K_0(Var_{\cc})$ to be the universal Euler characteristic.
After inverting the class $\mathbb{L}=[\mathbb{A}_{\cc}^{1}]$, we get the ring $\mathcal{M}_{\cc}$.

We introduce the equivariant Grothendieck group defined in \cite{DeLo1}.
Let $\mu_n$ be the cyclic group of order $n$, which can be taken as the algebraic variety
$Spec (\cc[x]/(x^n-1))$. Let $\mu_{md}\longrightarrow \mu_{n}$ be the map $x\mapsto x^{d}$. Then 
all the groups $\mu_{n}$ form a projective system. Let 
$$\underleftarrow{lim}_{n}\mu_{n}$$
be the direct limit.

Suppose that $X$ is a $S$-variety. The action $\mu_{n}\times X\longrightarrow X$ is called a $good$ 
action if  each orbit is contained in an affine subvariety of $X$.  A good $\hat{\mu}$-action on $X$ is an action of $\hat{\mu}$ which factors through a good $\mu_n$-action for some $n$.

The $equivariant ~Grothendieck~ group$ $K^{\hat{\mu}}_0(Var_{S})$ is defined as follows:
The generators are $S$-varieties $[X]$ with a good $\hat{\mu}$-action; and the relations are:
$[X,\hat{\mu}]=[Y,\hat{\mu}]$ if $X$ is isomorphic to $Y$ as $\hat{\mu}$-equivariant $S$-varieties,  
and $[X,\hat{\mu}]=[Y,\hat{\mu}]+[X\setminus Y, \hat{\mu}]$ if $Y$ is a Zariski closed subvariety
of $X$ with the $\hat{\mu}$-action induced from that on $X$,  if $V$ is an affine variety with a good 
$\hat{\mu}$-action, then $[X\times V,\hat{\mu}]=[X\times \mathbb{A}_{\cc}^{n},\hat{\mu}]$.  The group 
$K^{\hat{\mu}}_0(Var_{S})$ has a ring structure if we define the product as the fibre product with the good $\hat{\mu}$-action.  Still we let $\mathbb{L}$  represent the class $[S\times \mathbb{A}_{\cc}^{1},\hat{\mu}]$ and let $\mathcal{M}_{S}^{\hat{\mu}}=K^{\hat{\mu}}_0(Var_{S})[\mathbb{L}^{-1}]$ be the ring obtained from $K^{\hat{\mu}}_0(Var_{S})$ by inverting the class $\mathbb{L}$.

If $S=Spec(\cc)$, then we write $K^{\hat{\mu}}_0(Var_{S})$ as $K^{\hat{\mu}}_0(Var_{\cc})$, and $\mathcal{M}_{S}^{\hat{\mu}}$ as $\mathcal{M}_{\cc}^{\hat{\mu}}$.  Let  $s\in S$ be a geometric point. Then we have natural maps $K^{\hat{\mu}}_0(Var_{S})\longrightarrow K^{\hat{\mu}}_0(Var_{\cc})$ and $\mathcal{M}_{S}^{\hat{\mu}}\longrightarrow \mathcal{M}_{\cc}^{\hat{\mu}}$ given by the correspondence
$[X,\hat{\mu}]\mapsto [X_s,\hat{\mu}]$.

\subsection{The motivic Milnor fiber.}

Let $(L,\mu,\ae)$ be a cyclic $L_\infty$-algebra of dimension $3$.  Then we have a cyclic $L_\infty$-algebra structure $(H(L),\nu,\kappa)$ on the cohomology $H(L)$.  On $H^{1}(L)$ we have a formal series 
$$f: H^{1}(L)\longrightarrow \mathbb{C}.$$
Recall  that 
$$h: \YY\longrightarrow \MM$$
is the resolution of singularitiy of the formal scheme 
$\MM$.

Let $m_{I}=gcd(m_i)_{i\in I}$. Let $U$ be an affine Zariski open subset of $\YY$, such that, 
on $U$, $f\circ h=uv^{m_{I}}$, with $u$ a unit in $U$ and $v$ a morphism from 
$U$ to $\mathbb{A}_{\cc}^{1}$. The restriction of $E_{I}^{\circ}\cap U$, which we denote by
$\tilde{E}_{I}^{\circ}\cap U$, is defined by
$$\lbrace{(z,y)\in \mathbb{A}_{\cc}^{1}\times (E_{I}^{\circ}\cap U)| z^{m_{I}}=u^{-1}\rbrace}.$$
The $E_{I}^{\circ}$ can be covered by the open subsets $U$ of $Y$.  We can glue together all such 
constructions and get the Galois cover
$$\tilde{E}_{I}^{\circ}\longrightarrow E_{I}^{\circ}$$
with Galois group $\mu_{m_{I}}$.
Remember that $\hat{\mu}=\underleftarrow{lim} \mu_{n}$ is the direct limit of the groups
$\mu_{n}$. Then there is a natural $\hat{\mu}$ action on $\tilde{E}_{I}^{\circ}$.
Thus we get 
$[\tilde{E}_{I}^{\circ}]\in \mathcal{M}_{X_0}^{\hat{\mu}}$.

\begin{defn}
The motivic Milnor fiber of the cyclic $L_\infty$-algebra $L$ is defined as follows:
$$MF_0(L):=\sum_{\emptyset\neq I\subset A}(1-\mathbb{L})^{|I|-1}[\tilde{E}^{\circ}_{I}\cap h^{-1}(0)].$$
\end{defn}
It is clear that $MF_0(L)\in \mathcal{M}_{\cc}^{\hat{\mu}}.$

\begin{rmk}
Let $x\in X_0$ be a point and $f$ is convergent at $x$.  Then one can define 
 $$MF_x(L):=\sum_{\emptyset\neq I\subset A}(1-\mathbb{L})^{|I|-1}[\tilde{E}^{\circ}_{I}\cap h^{-1}(x)].$$
It is clear that $MF_x(L)\in \mathcal{M}_{\cc}^{\hat{\mu}}$.
\end{rmk}

\textbf{A constructible function.}
Taking Euler characteristic we get an integer value function:
\begin{equation}
\chi_{top}: X_0\longrightarrow \zz
\end{equation}
which is given by
$$x\mapsto \chi_{top}(MF_{x}(L)):=\sum_{\emptyset\neq I\subset A} \chi_{top}(\tilde{E}^{\circ}_{I}\cap h^{-1}(x))$$
for $x\in X_0$.
\begin{prop}
The function $\chi_{top}$ is a constructible function on $X_0$.
\end{prop}
\begin{pf}
This is easily from the fact that the function is given by taking Euler 
characteristic of the fiber of the 
morphism $\YY\rightarrow \MM$.
\end{pf}

\subsection{The analytic Milnor fiber.}

Let $R=\cc[[t]]$ be the formal power series ring, which is a discrete valuation ring. Let 
$K=\cc((t))$ be the quotient field. 
Let $\MM:=\widehat{H^{1}(L)}$ be the formal scheme.
In this section we recall the analytic Milnor fiber $\FF_{x}$ for any
$x$ in the special fiber of $\MM$ defined in \cite{NS}.

First note that the field $K$ is a non-archimedean field. And $K$ is also a discrete valuation ring with the 
standard valuation. The norm is given by 
$$|a|=c^{\val(a)}$$
for $a\in K$, where $c\in (0,1)$ is a fixed number.
Then   the norms of the elements of
$\cc$ are equal to $1$,  the formal series $f\in \cc[[x_1,\cdots,x_m]]$,  taken as a formal series in $K[[x_1,\cdots,x_m]]$, is convergent in the domain 
$U\subset (\aaa_{\cc}^{m})^{an}=\{|x_i|<1, 1\leq i\leq m\}$. 

Let 
$$\hat{f}_{x}: \MM\longrightarrow \cc[[t]]$$
be the morphism given by the formal series $f$, which is the formal completion of $f$ at $x$.
Let $\MM_0$ be the special fibre, and $\MM_{\eta}$ the generic fiber.  The generic fiber $\MM_{\eta}$ can be taken in the category of separated quasi-compact rigid $K$-varieties. Then from \cite{NS},
there is a canonical specialization morphism 
$$sp: \MM_{\eta}\longrightarrow \MM.$$
\begin{defn}
For any $x\in \MM_0$, the analytic Milnor fiber $\FF_{x}$ is defined as 
$\FF_{x}=sp^{-1}(\{x\})$.
\end{defn}

Let $\widehat{\MM\setminus \{x\}}$ is the formal completion of $\MM$ along 
$x$, then $\FF_{x}$ is isomorphic to the general fiber of the morphism $\widehat{\MM\setminus \{x\}}\rightarrow \MM$. Or as explained in \cite{KS2}, the analytic Milnor fibration is given by
\begin{multline*}
B(\epsilon)\cap f^{-1}(\eta)=\{x=(x_1,\cdots,x_m)\subset U | ~\mbox{max}_{i} |x_i|\leq \epsilon,
0<|f(x)|\leq\eta\} \\
\longrightarrow\{\omega\in(\aaa_{\cc}^{1})^{an}| 0<|\omega|\leq\eta\},
\end{multline*}
where $\eta\ll\epsilon\ll 0$. The analytic Milnor fiber is the generic fiber of such fibration.

From \cite{NS}, $MF(L)$ is actually the motivic Serre invariants and taking Euler characteristic
yields
\begin{equation}\label{analytic-motivic}
\chi_{et}(\FF_{x})=\chi_{top}(MF_{x}(L)),
\end{equation}
where $\chi_{et}(\FF_{x})$ is the Euler characteristic of the \'{e}tale cohomology of the analytic space
$\FF_{x}$ over the non-archimedean field $K$ in the sense of Berkovich \cite{Ber}.

\begin{rmk}
The homotopy type $\Delta(\FF_{x})$ of the analytic Milnor fiber $\FF_{x}$ was discussed by Nicaise \cite{Ni2}.   Getzler \cite{Getzler} defined a simplicial set $MC_{\bullet}(L)$ for any nilpotent $L_\infty$-algebra $L$, which is a Kan complex.  We conjecture that the homotopy type $MC_{\bullet}(L)$ of the 
cyclic $L_\infty$-algebra $L$ is equivalent to the homotopy type $\Delta(\FF_{x})$ of the analytic Milnor fiber $\FF_{x}$. We hope to return this later.

In general the superpotential function $f: H^1(L)\to \cc$ obtained  from a cyclic $L_\infty$-algebra 
is a formal power series.  In this paper we assume that $f$ is convergent at the point in $f^{-1}(0)$.
As in Kontsevich and Soibelman \cite{KS}, one can take the formal power series $f$ as a series over the non-archimedean field $K:=\cc((t))$ with a norm given by the standard valuation.  Then $f$ is convergent over a neighborhood of zero on the Berkovich analytic space $(\mathbb{A}^{m})^{an}$.  According to the theory of Berkovich \cite{Ber3} on vanishing cycles over non-archimedean field we may construct the vanishing cycle sheaf  of cyclic $L_\infty$-algebras using Berkovich theory, such that the stalk of the vanishing cycle sheaf is the Euler characteristic of the analytic Milnor fiber defined in this paper.  We leave this as a future work.
\end{rmk}

%%%%%%%%%%%%%%%%%%%%%%%%%%%%%%%%%%%%%%%

\section{The relation to Behrend function.}\label{behrend}

In this section we relate the motivic Milnor fiber to Behrend function defined in \cite{Be}.
\subsection{The Behrend constructible function.}

The cohomology $H^{i}(L)\neq 0$ only for $i=1,2$, where $H^{1}(L)$ classifies the deformation and $H^{2}(L)$ classifies the obstruction. Since $H^{1}(L)\cong H^{2}(L)^{\vee}$, from \cite{Be}, the obstruction theory associated to $L$ is symmetric.

There is a cyclic $L_\infty$-algebra structure $(H(L),\nu,\kappa)$ on the cohomology. 
The Maurer-Cartan space is defined by 
$$MC(L)=\{x\in H^{1}(L)| \sum_{n=2}^{\infty}\frac{(-1)^{n(n+1)}}{n!}\nu_{n}(x,\cdots,x)=0\}.$$

Let $X=MC(L)=Z(df)\subset H^{1}(L)$ be the Maurer-Cartan space, which is the germ moduli space of the deformation theory associated $L$. 
Suppose that $X$ is closed inside $M$, which means that $f$ is convergent at 
$x\in X$.
 As in \cite{Be}, let $C_{X/H^{1}(L)}$ be the normal cone, which belongs to the cotangent bundle $\Omega_{H^{1}(L)}$. Let
$$\pi: C_{X/H^{1}(L)}\longrightarrow X$$
be the projection. Let $C^{\prime}$ represent the irreducible components of $C_{X/H^{1}(L)}$.
Let the cycle $\Cc_{X/H^{1}(L)}$ on $X$ given by
$$\Cc_{X}:=\Cc_{X/H^{1}(L)}=\sum_{C^{\prime}}(-1)^{\mbox{dim} \pi(C^{\prime})}\mbox{mult}(C^{\prime})\pi(C^{\prime}).$$
Here $\mbox{mult}(C^{\prime})$ is the length of $C_{X/H^{1}(L)}$ at the generic point of $C^{\prime}$.

The Behrend canonical constructible function $\nu_{X}$ is defined to be the Euler obstruction function 
$\mbox{Eu}(\Cc_{X/H^{1}(L)})$. The Euler obstruction $Eu_{\Cc_{X}}$ is a constructible function 
that assigns the value 
$$Eu_{\Cc_{X}}(P)=\int_{\mu^{-1}(P)}c_{*}(\tilde{T})\cap s(\mu^{-1}(P),\tilde{\Cc_{X}})$$
for any point $P\in X$, where $\mu: \tilde{\Cc_{X}}\longrightarrow \Cc_{X}$ is the Nash blowup and $\tilde{T}$ is the universal tangent bundle. For details, see \cite{Be}.

The local moduli space $X$ is determined by the cyclic $L_\infty$-algebra $L$, we have:
\begin{prop}
The germ moduli space $X$ admits a symmetric obstruction theory in the sense of \cite{Be}.
\end{prop}
\begin{pf}
The deformation is classified by $H^{1}(L)$, the first cohomology, and the obstruction is given by the 
second cohomology $H^{2}(L)$. By cyclic property, $H^{1}(L)\cong (H^{2}(L))^{\vee}$, and the deformation is isomorphic to the dual of the obstruction.  Hence the obstruction theory is symmetric.
\end{pf}

\textbf{Characteristic cycles.}

Let 
$$L:  Z_{*}(X)\longrightarrow \mathcal{L}_{X}(\Omega_{H^{1}(L)})$$
be the map sending the integral value  cycles to conic Lagrangian cycles.  Then
from \cite{Be},  
$$L(\Cc_{X/H^{1}(L)})=C_{X/M},$$ is the normal cone.

Let 
$$L:  \mbox{Con}_{*}(X)\longrightarrow \mathcal{L}_{X}(\Omega_{H^{1}(L)})$$
be the map sending the  constructible functions to conic lagrangian cycles.  Then 
\begin{equation}\label{nu-c}
L(\nu_{X})=C_{X/M}.
\end{equation}

Let $\jeden_{V}$ be the canonical constructible function for a subvariety $V\in X$. Then 
if $V=\sum_{i}a_iY_i$, for subvarieties $Y_i\subset V$, we have 
$\jeden_{V}=\sum_{i}a_i\mbox{Eu}(Y_i)$.

Recall that there is  a formal series 
$$f: H^{1}(L)\longrightarrow \mathbb{C}$$
given by the cyclic relations $\nu$ and the binear form $\kappa$. 
If the series $f$ is a holomorphic function on $M$, then 
$$\nu_{X}(0)=(-1)^{\dim(H^{1}(L))}(1-\chi(F_{0}))$$
where $F_{0}$ is the Milnor fibre of the function of $f$ at the point $0\in H^{1}(L)$.

%%%--------------------------------------------------------------
\subsection{The formal blow-up.}

Let $M:=H^{1}(L)$. And let $I\subset\O_{M}$ be the ideal of the 
singular subvariety of $M$ determined by the differential $df$.

Then from \cite{Tem} and \cite{NS}, we have 
\begin{prop}
The resolution $h: \YY\rightarrow \MM$ can be obtained by blow-up
along the vanishing locus $\mathbb{V}(I)$ and the completion of the 
blow-up is naturally isomorphic to the formal blow-up. i.e.
$$\widehat{\mbox{Bl}_{\mathbb{V}(I)}M}\cong \mbox{Bl}_{\mathbb{V}(df)}\mathfrak{M}.$$
\end{prop}

Then we have:

\begin{lem}\label{resolution-conormal}
The  resolution $\mbox{Bl}_{\mathbb{V}(I)}M$ of $(M, X_0)$ is isomorphic to the projectivization of the 
conormal bundle $T_{f}^{\vee}M$ inside $\Omega_{M}$, hence
$$E_i=\mathbb{P}(T_{h(E_i)}^{\vee}M).$$
\end{lem}
\begin{pf}
The closed subscheme $X\subset M$ is the singular locus of the scheme $M$.
From \cite{Tem}, the  resolution $h: \YY\longrightarrow \MM$ can be obtained from  
blowups along the singular locus $X$.
\end{pf}

The formal scheme $\MM\rightarrow \cc[[t]]$ is defined over the formal power series ring 
$\cc[[t]]$. 
We take the  formal series $f$ in the non-archimedean field $K$, hence define a
global section of $\Omega_{\MM_{\eta}/K}^{0}(\MM_{\eta})$ for the generic fiber $\MM_{\eta}$. 
In the 
resolution $h$, $m_i$ is the order of $h^{*}f$ along the exceptional divisor
$E_i$, which can defined as the length of the $\O_{\mathfrak{M},\xi}$-module
$\O_{\mathfrak{M},\xi}/f\O_{\mathfrak{M},\xi}$ at a generic point 
$\xi\in E_i$ in the formal scheme $\mathfrak{M}$.  
Similarly let 
$n_i$ be the order of $(f)$ along $h(E_i)=C_i$ and $p_i$ the order 
of the Jacobi ideal $(df)$ along $E_i$. 

\begin{prop}\label{orders} 
$$m_i=n_i+p_i.$$ 
\end{prop}
\begin{pf}
The Jacobi ideal $(df)$ is given by the cokernel of the induced morphism
$$h^{*}\Omega_{\MM}\longrightarrow \Omega_{\YY}.$$
Let $\gamma(t)$ be an analytic curve in the exceptional divisor 
$E$ over the $\cc((t))$. We compute that
$$\ord_{0}(f(\gamma(t)))=\ord_{0}(x(t))+\ord_{0}\left(\frac{\partial f}{\partial x_1}(x_1(t)),\cdots, \frac{\partial f}{\partial x_m}(x_m(t))\right)$$
implies the proposition. 
Then this is a routine computation as in \cite{PP}.
\end{pf}

\begin{rmk}
This proposition can be taken as a generalization of Proposition 2.2  in \cite{PP} for a regular function 
$f$ to the formal power series inside the category of rigid spaces over non-archimedean field.
\end{rmk}

%%%---------------------------------------------------------------------------
\subsection{Relation to Behrend function.}

Recall that the Behrend function 
$$\nu_{X}: X\rightarrow \zz$$
satisfies
$$L(\nu_{X})=C_{X/M}\subset \Omega_{M},$$ 
as in (\ref{nu-c}).

\begin{lem}\label{C=E}
The cycle $C_{X/M}$ is equivalent to the cycle $\sum_{i}p_i E_i$ in $\Omega_{M}$.
\end{lem}
\begin{pf}
Let $h: \YY\longrightarrow \MM$ be the blow-up along the singular subscheme $X\subset M$.
Let $I\subset \O_{M}$ be the ideal of $X$.
From the definition of blow-up, the projectivization of the normal cone 
$\mathbb{P}(C_{X/M})=\mathbb{P}(\bigoplus_{n} I^{n}/I^{n+1})$ is equal to 
the exceptional divisor $E$.

By Lemma \ref{resolution-conormal}, the resolution $h$ is isomorphic to the projectivization of the   conormal bundle, and 
hence $\pp(C_{X/M})=\pp(\sum_{i}p_i E_i)=E$.
\end{pf}

\begin{prop}\label{11}
$$L(\nu_{X})=(-1)^{m-1}(L(\chi_{top})-L(\jeden_{X_0})).$$
\end{prop}
\begin{pf}
From $E_i=\mathbb{P}(T_{C_i}^{\vee}M)$,  and 
$$L(\jeden_{X_0})=(-1)^{m-1}\sum_{i}n_i E_i , ~~L(\chi_{top})=(-1)^{m-1}\sum_{i}m_i E_i$$
and Proposition \ref{orders}, Lemma \ref{C=E},  the result follows.
\end{pf}

Then from the correspondence of the constructible functions and 
Lagrange cycles one has:
\begin{cor}
$$\nu_{X}=(-1)^{m-1}(\chi_{top}-\jeden_{X_0}).$$
\end{cor}

Over the origin $0\in M$, we have
$$
\nu_{X}(0)=(-1)^{m}\left(1-\sum_{i}m_i\chi(E_i\cap h^{-1}(0))\right).
$$
From (\ref{analytic-motivic}), we have
\begin{cor}
We have
\begin{equation}\label{formal-neighborhood}
\nu_{X}(0)=(-1)^{m}\left(1-\chi_{et}(\FF_{0})\right),
\end{equation}
where $\FF_0$ is the analytic Milnor fiber over $0\in \MM$ in the sense of 
\cite{NS}.
\end{cor}
\begin{rmk}
The formula (\ref{formal-neighborhood}) means that the value of the Behrend function 
$\nu_{X}$ is given by the analytic Milnor fiber in the formal setting, hence is independent  of the formal neighborhood of $0$.
\end{rmk}

\subsection{Regular function.}
In this section we assume that the series $f$
is a regular function. 
Then there is a morphism $f: M\rightarrow  Spec(\cc[t])$
and the formal variety $\MM$ is the $t$-adic completion of $f$.
\begin{rmk}
If the cyclic $L_\infty$-algebra $L$ is from a Schur object in the 
derived category of coherent sheaves on Calabi-Yau 3-folds,  the formal potential $f$ has been already proved by Behrend-Getzler \cite{BG} to be convergent.
\end{rmk}

We still let $X=\mathbb{V}(df)$.  For the origin  $0\in M$,
let $F_0$ be the Milnor fiber 
 at the point $0$ in the sense of \cite{Mil}.  Then from A'campo's \cite{AC} formula, 
 which states that the Euler characteristic of the exceptional divisor is the Euler characteristic 
 of the Milnor fiber, or the result  
 $\chi_{et}(\FF_{0})=\chi_{top}(F_{0})$ of \cite{NS},
 the following result is obvious.
$$\nu_{X}(0)=(-1)^{m}(1-\chi_{top}(F_{0})),$$
where $\chi_{top}(F_0)$ is the topological Euler characteristic of the Milnor fiber.

%%%----------
\textbf{Relation to Donaldson-Thomas invariants.} 

Let $L$ be a Donaldson-Thomas type $\dg$ Lie algebra, i.e. $L$ is a cyclic $\dg$ Lie algebra 
of dimension three such that $H^{i}(L)=0$ only except $i=1,2$.  
Let $0\in M$ be the origin.
Then from \cite{Be} and \cite{BG},  $\nu_{X}(0)$ is the pointed Donaldson-Thomas invariant.
Then 
\begin{thm}
  Let $0\in H^{1}(L)$ be the origin, then the value
$\nu_{X}(0)$ is the pointed Donaldson-thomas invariant and is equal to the Euler characteristic 
of the motivic Milnor fiber of the cyclic dg Lie algebra $L$. $\square$
\end{thm}

%%%%%%%%%%%%%%%%%%%%%%%%%%%%%%
%%%----------------------------------------------------------------------

\subsection*{}

% ------------------------------------------------------------------------
\end{document}